\documentclass[12pt]{amsart}

\usepackage{amssymb,amsmath,amsthm}
\usepackage{graphicx}
\usepackage[mathscr]{eucal}

\date{October 29, 2008}
\usepackage{geometry}
\usepackage{graphicx}
\DeclareMathOperator{\bd}{bd}
\DeclareMathOperator{\Cat}{Cat}
\DeclareMathOperator{\R}{R}
\DeclareMathOperator{\dia}{dia}
\DeclareMathOperator{\starr}{star}
\newcommand{\cat}{\Cat(k)}
\newcommand{\catt}{\Cat(0)}
\newcommand{\sett}[2]{ \{#1, \dots ,  #2 \}  }  
\newcommand{\model}{M^2_k}
\newcommand{\bld}{\mathbb{R}}
\newcommand{\fn}[3]{#1:#2 \rightarrow{#3}}
\newcommand{\tri}[3]{\bigtriangleup(#1,#2,#3)}
\newcommand {\tribar}[3]{ \bigtriangleup ( \bar{#1},\bar{#2},\bar{#3} ) }
\newcommand{\inff}[4]{\underset{#1}{ \inf} \sup \{#2 \mid #3 \in #4 \}}
\newcommand{\limm}[2]{\underset{l \to \infty}{\lim }\sup \{u(p)d(#1,p) \mid p \in #2    \} }     
\newcommand{\m}[2]{m_k(#1,#2)}
\newcommand{\disttt}[2]{\underset{#1}{ \inf} \sup \{u(p) d(#2,p) \mid p \in P \}}
\newcommand{\dis}[2]{\underset{#1}{ \inf} \sup \{ {#2} \mid p \in P \}}

\newcommand{\distt}[4]{\sup \{u(#1) d(#2,#3) \mid #4  \}}
\newcommand{\fdelta}{f_{\delta}}
\newcommand{\ball}[2]{B(#1,#2)}
\newcommand{\oball}[2]{\overline{B(#1,#2)}}

\newcommand{\script}[1]{\mathscr{#1}}
\newcommand{\midd}[1]{\overline{m_{#1 ll'} }}

\newtheorem{lemma}{Lemma}
\newtheorem{theorem}{Theorem}
\newtheorem{corollary}{Corollary}
\newtheorem{definition}{Definition}
\newtheorem{proposition}{Proposition}
\newtheorem{example}{Example}

\title{A Generalization of the Circumcenter of a Set}

\author{Jack E. Girolo}
\address{Department of Mathematics, California Polytechnic State University, 1 Grand Avenue, San Luis Obispo, California 93407}
\email{jgirolo@calpoly.edu}


\begin{document}

\maketitle

\begin{abstract}
Let $(X,d)$ be a $\cat$ space and  $P$  a bounded subset of $X$.
If $k>0$ then it is required that the diameter of $P$ be less than $\frac{\pi}{4\sqrt{k}}$. Let $u:P\rightarrow \mathbb{R}$ be a bounded non-negative function from $P$ to $\mathbb{R}$. The existence of a unique point in $X$ called the barycenter of $P$ relative to $u$ is established. When $u\equiv 1$, the barycenter is simply the circumcenter of $P$. The barycenter has a number of properties including a scaling, continuity and limit property. Under suitable conditions, the barycenter is a fixed point of an isometry or group of isometries. Barycenters are used to show that a complete $\cat$ space $X$ is an absolute retract if $k \leqq 0$, and an absolute neighborhood retract if $X$ is complete and of curvature $\leq k$.
\end{abstract}

\section {Introduction. }\label{S:intro} 
 Given a bounded set $P$ in $\mathbb R^n$, the classical Jung Theorem \cite{Jung} states that there is a unique
 smallest closed ball $B(q_1,r_1)$ containing $P$, where
 
  $\displaystyle{radius(B) \leqq \sqrt{\frac{n}{2(n+1)}}diameter(P)}$.  
   \newline Furthermore, equality holds if and only if $\overline{P}$, the closure of $P$, contains a regular n-simplex with diameter equal to the diameter of P.  The uniqueness statement above can be formulated as follows.
 There exists a unique $q_1$ in $\mathbb{R}^n$ and $r_1$ in $\mathbb{R}$ such that

$\displaystyle{ r_1= \sup\{d(q_1,p) \mid p \in P\}=\inff{x}{d(x,p)}{p}{P} }$. \newline
The point $q_1$ is called the circumcenter of $P$ and the number $r_1$ is called the circumradius of $P$.

 Let $u:P\rightarrow \mathbb{R}$ be a bounded non-negative function from $P$ to $\mathbb{R}$. The function $u$ is not required to be continuous. We prove that there exists a unique $q_u$ in $\mathbb{R}^n$ and $r_u$ in $\mathbb{R}$ such that 
 
 $\displaystyle{ r_u= \sup\{u(p)d(q_u,p) \mid p \in P\}=\inff{x}{u(p)d(x,p)}{p}{P} }$. \newline
 The point $q_u$ is called the barycenter of $P$ relative to $u$ and $r_u$ is called the baryradius of $P$ relative to $u$. It will be shown that $d(q_1,q_u) \leqq r_1$.
 
 Jung's Theorem has been extended to the $\cat$ setting by Lang and Schroeder \cite{ls}, and it is in this setting that the concept of the barycenter and baryradius are developed. When $u \equiv 1$, the barycenter and baryradius are simply the circumcenter and circumradius of $P$. The barycenter has a number of properties which include scaling, continuity and limit properties. Under suitable conditions, the barycenter is a fixed point of an isometry or group of isometries.
 
 Convex sums in linear spaces, when combined with partitons of unity, are an important tool in constructing the extension of continuous functions (\cite{gd}, P.163-164). We show that the barycenter (of finite sets) plays an analogous role in the $\cat$ setting. In particular we show that  $\cat$ spaces, $k\le 0$, are absolute retracts, $\cat$
spaces, $k > 0$, are absolute neighborhood retracts, and spaces of curvature less than or equal to $k$ are absolute neighborhood retracts. For $X$ a $\catt$ space, we state a result similar to Michael's Selection Theorem and, given a group of isometries $\Gamma$ with bounded orbits, we show that  a bounded function $\fn{u}{X}{\bld}$ that is positive valued and uniformly continuous induces a continuous function on the orbit space $\Gamma/X$.

 \section{ Definitions and Basic Results.}\label{S:basic} The spaces that we consider had their origin in the work of 
A.~D.~Alexandrov~\cite{alex51}, in which he extended the fundamental concept of curvature in Riemannian manifolds  
 to a subclass of geodesic metric spaces called $\cat$ spaces. In this paper all $\cat$ spaces are assumed to be complete. This he accomplished by comparing triangles in geodesic metric spaces to triangles in a specific class of Riemannian manifolds of constant curvature, one with curvature $k$ for each real number $k$. The expression $\cat$ was given by Gromov~\cite{gro87} in recognition of the fundamental work done in this subject by Alexandrov, E.~Cartan, and A.~Toponogov.

 The recent book by Martin R. Bridson and 
Andr\'{e} Haefliger~\cite{bh} is a fundamental source for this subject. It is from their book that we 
take the basic definitions and results pertinent to this paper.

By a geodesic metric space we mean a metric space $(X,d)$ such that for any two points $x$ and $y$
in $X$ there exists a path, a continuous function, $\sigma:[0,l] \rightarrow X$ from the closed interval
$[0,l]$ to $X$ such that $\sigma(0)=x$, $\sigma(l)=y$ and $d(\sigma(t),\sigma(t'))=\left|t-t'\right|$ for $t,t' \in [0,l]$. In particular $d(x,y)=l$ In this paper  $l$ is finite. Such a path will be called a geodesic segment.
 
 We introduce three different classes of model spaces upon which the definition  of a $\cat$ space is based. For $\sett{x_1}{x_n}$ and $y=\sett{y_1}{y_n}$ in $\mathbb R^n$, define 
 $x\vert_E y= \sum_{i=1}^{n}x_i y_i $, the Euclidean scalar product. Set $(E^n,d_E)$ to be $ \mathbb R^n$ with the Euclidean norm $ \left\| x \right \| =\sqrt {x|_Ex}$ and metric $d_E(x,y)= \left\| x-y \right \| $.
Next let $S^n = \{ x \in E^{n+1} \mid \left\| x \right \| =1 \}$ with metric $d_S$ given by  $ \cos(d_S(x,y))=x|_Ey$. Finally, let $x=\sett{x_1}{x_{n+1}}$ and $ \sett{y_1}{y_{n+1}}$ be in $\mathbb R^{n+1}$ and let $x|_Hy$ be the symmetric bilinear form given by 
$x|_Hy=\sum_{i=1}^{n}x_iy_i-x_{n+1}y_{n+1}$. Set $H^n = \{ x \in \mathbb R^{n+1}\mid x|_Hx=-1
\text{ and } x_{n+1}>0 \}$. Then $n$ dimensional hyperbolic space is given by $(H^n,d_H)$ where 
$\cosh(d|_H(x,y)))=x|_Hy$.
\paragraph{The Model Spaces $M^n_k$}  For $k>0$, $M_k^n=(S^n,d^k_S)$ with
$\displaystyle{d^k_S(x,y)=\frac{1}{\sqrt{k}}d_S(x,y)}$, for $k=0$, $M^n_0=(E^n,d_E)$ and
for $k<0$, $M^n_k=(H^n,d^k_H)$ with
$\displaystyle{d^k_H(x,y)=\frac{1}{\sqrt{-k}}d_H (x,y)}$.

In that which follows, the metric in every metric space will be denoted by the expression $d(x,y)$.
In particular the subscript and superscript notation above will be omitted . Following standard notation as found in \cite{bh}, if $k>0$ set $D_k= \pi/\sqrt{k} $, the diameter of the model space $M^n_k$. For $k\le 0$ the model space $M^n_k$ has an unbounded metric and we set $D_k=\infty$. In a geodesic
metric space $(X,d)$, a subset $C$ is convex provided that for any two points $x$ and $y$ in $C$ a
geodesic segment joining $x$ to $y$ lies completely in $C$. In the following two propositions we gather some of the basic properties of the model spaces.

\begin{proposition}

The space $M^n_k$ is a geodesic metric space. If $k\leq0$ then $M^n_k$ is 
uniquely geodesic and all balls in $M^n_k$ are convex. If $k>0$ then there is a unique geodesic
segment joining $x$ to $y$ in $M^n_k$ if and only if $d(x,y)<D_k$. If $k>0$ then closed balls in $M^n_k$ of radius less than $ D_k/2$ are convex.
\end{proposition} 

A geodesic triangle in a geodesic metric space $(X,d)$ consists of three points $x,y$ and $z$, and 
three geodesic segments $xy$, $yz$ and $xz$  joining each pair of the three points. Recall that the length of $xy$, for example, equals the distance from $x$ to $y$. Such a triangle will be denoted by $\bigtriangleup(xy,yz,xz)$ or more simply, but less precisely, by $\bigtriangleup(x,y,z)$. Typical points in a model space will be denoted by $\bar{x},\bar{y}$ and $\bar{z}$, and the corresponding geodesic triangle by $ \bigtriangleup ( \bar{x},\bar{y},\bar{z} ) $. The law of cosines for each of the model spaces follows. 

\begin{proposition} Let  $\bigtriangleup$ be a triangle in $\model$ with sides of length $a,b$ and $c$, and let $\alpha$ be the angle opposite the side of length $a$. Then 
\begin{itemize}
 \item for $k>0$, $\cos(\sqrt{k}a)= \cos(\sqrt{k}b)\cos(\sqrt{k}c)+\sin(\sqrt{k}b)\sin(\sqrt{k}c)\cos(\alpha)$,
\item for $k=0$, $a^2=b^2 +c^2-2bc\cos(\alpha)$, and 
 \item for $k<0$, $\cosh(\sqrt{-k}a)= \cosh(\sqrt{-k}b)\cosh(\sqrt{-k}c)-\sinh(\sqrt{-k}b)\sinh(\sqrt{-k}c)\cos(\alpha)$.
\end{itemize}

\end{proposition}

While we have defined model spaces in  each dimension we will only use those of dimension $2$
in the definition of $\cat$ spaces. Let $\tri{x}{y}{z}$ be a geodesic triangle in a geodesic metric space $(X,d)$. 
Then a geodesic triangle $\tribar{x}{y}{z}$ in  $\model $ is a comparison triangle of $\tri{x}{y}{z}$ provided
$d(\bar{x},\bar{y})=d(x,y),d(\bar{y},\bar{z})=d(y,z)$ and $d(\bar{x},\bar{z})=d(x,z)$. Suppose that $x'$ is an element of $xy$. Then $\bar{x'}$ in $\tribar{x}{y}{z}$ is the comparison point for $x'$ provided that
$\bar{x'}$ is in $\bar{x}\bar{y}$ and $d(x',x)=d(\bar{x'},\bar{x})$. Comparison points for the other sides
of the triangle are similarly defined.

\begin{definition}\label{D:cat}   Let $(X,d)$ be a geometric metric space. Then $X$ is a $\cat$ space  provided that for every geodesic triangle $\tri{x}{y}{z}$ in $X$ with perimeter less than $2D_k$, there exists a comparison triangle $\tribar{x}{y}{z}$ in $\model$ such that for any two points $x'$ and $y'$ in $\tri{x}{y}{z}$ and comparison points $\bar{x'}$ and $\bar{y'}$ in $\tribar{x}{y}{z}$ the inequality $d(x',y')\leq d(\bar{x'},\bar{y'})$ is true. This inequality is called the $\cat$ inequality.
\end {definition}

Riemannian manifolds of bounded sectional curvature are $\cat$ spaces (\cite{bh},P.173). It follows that for each $k$ the model space $M^k_n$ is a $\cat$ space. A wide variety
of spaces are $\cat$ as illustrated by the following examples (\cite{bh}, P.167-168). 

\begin{example}\label{E:space}

\begin{enumerate}
\item A normed linear space is a $\cat$ if and only if it is
a pre-hilbert space. 
\item A subset of a $\cat$ space is $\cat$ if and only if it is convex.
\item In a $\cat$ space, a ball of diameter less than $D_k$ is a $\cat$ space.
\item A metric simplicial graph is a $\cat$ 
space if and only if it does not contain an essential
loop of length less than $2D_k$. 
\item An $R$-tree is $\cat$ for every $k$.
\end{enumerate}
 \end{example}
 
 There are a number of equivalent definitions of $\cat$ spaces, some of which are conceptually less intricate, and a number of properties that we do not state in this paper. We refer the reader to \cite{bh} and also \cite{bbi}.
 
 \section{Barycenters.}\label{S:coor}

Given a complete $\cat$ space $(X,d)$, a subset $P$ of $X$ of diameter that is finite and less than $D_k/4$, and a bounded and non-negative function $u:P\rightarrow \mathbb R$, we establish the existence of unique point that can be considered as the barycenter of $P$ corresponding to the function $u$. Furthermore, we note that no continuity condition is imposed upon the function $u$. We will see that this barycenter has a number of useful properties.  

\begin{lemma}\label{L:inequal1} Let $xyz$ be a geodesic segment in $\model$
such that $r=d(x,y)<d(x,z)\leq D_k/2$. Let $w$ be in $\overline{B(x,r)}$. Then $ d(y,w)<d(z,w)$.
\end{lemma}

\begin{center}
\includegraphics[width=2in]{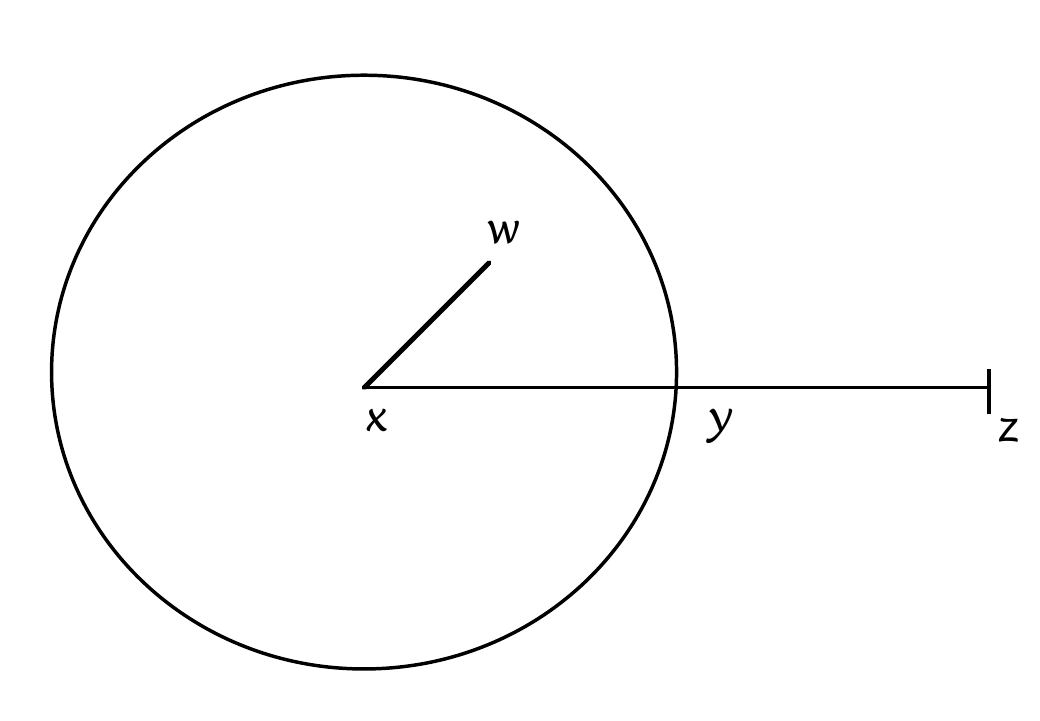}\\
Figure 1.
\end{center}

\begin{proof} For each $k$ an application of the law of cosines shows that $\angle xyw$ and  $\angle xzw$ are acute. It follows that $\angle wyz$ is obtuse. Consider $\tri{w}{y}{z}$. Then $d(y,w)<d(z,w)$ as the greater side lies opposite the greater angle. \end{proof}

Let $x$ be a point in $\model$. Let $r$ and $s$ be any two positive numbers satisfying $r<s< D_k/2$. We note that when $k\leq 0$ no upper bound is imposed upon $s$. Define $ W=\overline{B(x,r)}$, $Y=\{y \mid d(y,x)=r\}$ and $Z=\{z \mid d(x,z)=s\} $.  Let $\pi$ be the projection along geodesic segments of $Z$ onto $Y$, and define $f: W \times Z \rightarrow\mathbb{R}$ by $f(w,z)=d(w,z)-d(w,\pi(z))$.

\begin{lemma}\label{L:inequal2} The function $f :W \times Z \rightarrow{\mathbb{R}}$, as defined above, has a minimum value denoted by $m_k(r,s)$, and
$m_k(r,s)>0$.
\end{lemma}
\begin{proof} The set $W \times Z$ is compact. Since $f $ is continuous, it has a minimum value. So there is a $w$ in $W$ and a $z$ in $Z$ 
such that $m_k(r,s)=f(w,z)=d(w,z)-d(w,\pi(z))$. By Lemma \ref{L:inequal1} $m_k(r,s)>0$.
\end{proof}

\begin{definition}{\label{D:disttt}} Let $X$ be a $\cat$ space, and $P$ be a subspace of  $X$ with diameter that is finite and less than $D_k/4$. Let $u:P\rightarrow \mathbb R$ be a bounded and non-negative real valued function from $P$ to $\mathbb R. $\footnote{We assume that $u$ is not identically $0$.} Then define $r_u=\disttt{x \in X}{x}$.

\end{definition}

\begin{lemma}\label{L:rupos} Let $X$, $P$, $u$, and $r_u$ be as in Definition
\ref{D:disttt}. If there exists two distinct points $p_1$ and $p_2$ in $P$ such that $u(p_1)>0$ and $u(p_2)>0$, then $r_u>0$. 
\end{lemma}

\begin{proof}  Let $x$ be an arbitrary element of $X$. Without loss of generality, we may assume that $d(x,p_1) \geq d(p_1,p_2)/2$. Then $\distt{p}{x}{p}{p \in P} \geq \min\{ u(p_1),u(p_2) \} d(p_1,p_2)/2$. This inequality being true for each $x$ in $X$ implies that $r_u>0$. 

\end{proof} 

Let $x_l$ be a sequence in $X$ with the property that $\underset{l \to \infty} \lim \distt{p}{x_l}{p}{p \in P}=r_u$. On account of the definition of $r_u$, such a sequence exists. We estimate $d(x_l,p)$ for each $p \in P$ and for $l$ sufficiently large. In that which follows, $\dia(P)$ will denote the diameter of $P$.

\begin{lemma}\label{L:close} Let $X$, $P$, $u$, and $r_u$ be as in Definition
\ref{D:disttt}. Let $x_l$ be a sequence in $X$ with the property that $\underset{l \to \infty} \lim \distt{p}{x_l}{p}{p \in P}=r_u$. For every $\epsilon>0$ there exists a number
$N>0$ and a $p'$ in $P$ such that if $l>N$ then $d(x_l,p')< \dia(P) + \epsilon$ and, for every $p \in P$, $d(x_l,p)<2\dia(P)+\epsilon$.
\end{lemma}

\begin{proof} Let $\epsilon>0$ and set $b=\sup\{ u(p) \mid p \in P \}$. Since $u$ is a bounded function, $b<\infty$. Set $\dia(P)=r$. Pick $p$ in $P$. Then $r_u \leq \distt{p'}{p'}{p}{p'\in P} \leq br$. Since $b=\sup\{ u(p) \mid p \in P \}$ there exists a $p'$ in $P$ such that $u(p')(r+\epsilon)>br\geq r_u$. If for each $N>0$ there were to exist an $l>N$ such that $d(x_l,p')>r+\epsilon$, then $\distt{p}{x_l}{p}{p\in P}> u(p')(r+\epsilon)>r_u$. This is contrary to the fact that  $\underset{l \to \infty} \lim \distt{p}{x_l}{p}{p \in P}=r_u$. Therefore there exists an $N>0$ such that if $l>N$ then $d(x_l,p')<r+\epsilon=\dia(P)+\epsilon$. An application of the triangle inequality completes the proof.
\end{proof}

\begin{lemma}\label{L:P'} Let $X$, $P$, $u$, and $r_u$ be as in Definition
\ref{D:disttt}. Suppose that $z$ is in $X$, $r$ is positive and $P$ is a subset of $\oball{z}{r}$. Then there exists a subset $P'$ of $P$ and a number $a>0$ with the following properties: $a=\inf\{u(p) \mid p \in P' \}$ and for $x$ in $\oball{z}{r}$, $\distt{p}{x}{p}{p \in P'}=\distt{p}{x}{p}{p \in P}\geq r_u$. Furthermore, for $x$ in $\oball{z}{r}$ and $p \in P- P'$, $u(p)d(x,p) < r_u$.  
\end{lemma} 
\begin{proof} Choose $a$ so that $0<a<r_u/(2r+1)$. Set $P''=\{ p \in P \mid u(p) \leq a \}$ and $P'=P - P''$. Then for each $x$ in $\oball{z}{r}$ and $p$ in $P''$, $u(p)d(x,p)<\frac{r_u}{2r+1}2r<r_u$. So for $x$ in $\oball{z}{r}$,
$\distt{x}{x}{p}{p \in P'}=\distt{x}{x}{p}{p \in P} \geq r_u$.
\end{proof}

\begin{theorem}\label{T:center} Let $(X,d)$ be a complete $\cat$ space, and let $P$ be a subset of $X$ with diameter less than $D_k/4$. Let $u:P\rightarrow \mathbb R$ be a bounded non-negative function defined on $P$. Then there is a unique $q_u$ in $X$ such that $\distt{p}{q_u}{p}{p \in P}=\disttt{x}{x}$.
\end{theorem}

\begin{proof} Let $d_k$ be be a real number such that $\dia(P)<d_k<D_k/4$. Let $\{x_l\}$ be a sequence in $X$ such that $\limm{x_l}{P} =r_u$. We show that this sequence is Cauchy. 

By Lemma \ref{L:close} there exists an $N>0$ such that if $l>N$, then, for every $p$ in $P$,  $d(x_l,p)< 2d_k$.

Choose a $p'$ in $  P $. By Lemma \ref{L:P'} there exists an  $a>0$ and a subset $P'$ of $P$ such that $a=\inf\{ u(p) \mid p \in P'\}$ and for every $x$ in $B(p',2d_k)$, and in particular every $l>N$, $\distt{p}{x_l}{p}{p\in P'}=\distt{p}{x_l}{p}{p\in P}$. Therefore, $\limm{x_l}{P'} =r_u$. So it is sufficient to use $P'$  rather than $P$ to verify that the sequence $\{x_l\}$ is Cauchy. Toward that end, let $\epsilon>0$ and pick a $\bar{c}$ in $\model$.

For $k>0$ set $s=D_k/2-2d_k$ and choose $\delta>0$ so that 

1) $\delta/a < s/4$, and

2) if  $r_u/u(p)< 2d_k +s/4$, then a geodesic segment in the annulus

$\displaystyle{T=\{ \bar{x} \in \model \mid r_u/u(p)-\delta/a \leq d(\bar{c},\bar{x}) \leq r_u/u(p) +\delta/a \}}$ has length less than $\epsilon/2$ .
Applying the definition of $d_k$ and $s$, we see that

 $\displaystyle{2d_k+s/4=2d_k+\frac{D_k/2-2d_k}{4}=3d_k/2+D_k/8<D_k/2}$. So we can make the above claim on the length of geodesics in $T$.

For the case when $k\leq 0$, choose $\delta>0$ so that a geodesic segment in the annulus $\displaystyle{T=\{ \bar{x} \in \model \mid r_u/u(p)-\delta/a \leq d(\bar{c},\bar{x}) \leq r_u/u(p) +\delta/a \}}$ has length less than $\epsilon/2 $. We may do so since $\{ r_u/u(p) \mid p \in P' \}$ is bounded above.  

From the definition of $r_u$ and the definition of the sequence $\{x_l\}$ there exists an $M>N$ such that if $l$ and $l'$ are greater than $M$ then
 
 3) $r_u \leq \distt{p}{x_l}{p}{p \in P'}< r_u+\delta$ and $r_u \leq \distt{p}{x_{l'}}{p}{p \in P'}< r_u+\delta. $ 

Suppose that $l$ and $l'$ are greater than $M$. Let $m$ be the midpoint of the geodesic segment $x_lx_l'$. For $p' \in P'$ consider the geodesic triangle $\tri{p'}{x_l}{x_l'}$. On account of the definition of $d_k$ and the fact that both $l$ and $l'$ are greater than $N$, we have $d(x_l,p')<2d_k<D_k/2$ and $d(x_l',p')<2d_k<D_k/2$. By the triangle inequality, $d(x_l,x_l')<D_k$. So the perimeter of $\tri{p'}{x_l}{x_l'}$ is less than $2D_k$. Thus there is a comparison triangle $\tribar{c}{x_l}{x_l'}$ in $\model$. 

Let $\midd{p'}$ be the midpoint of the geodesic segment $\bar{x_l}\bar{x_l'}$. Suppose that for each $p' \in P'$, $u(p')d(\bar{c},\midd{p'})<r_u-\delta$. Then by the $\Cat(k)$ inequality $u(p')d(p',m)<r_u-\delta$ for each $p' \in P'$. It follows that $\distt{p'}{p'}{m}{p'\in P'} <r_u$. This is contrary to the definition of $r_u$. So there must exist a $p$ in $P'$ such that when one considers the comparison triangle $\tribar{c}{x_l}{x_l'}$, $u(p)d(\bar{c},\midd{p})>r_u-\delta$. For that which follows set $\midd{p}=\bar{m}$. From the inequality above $d(\bar{c},\bar{m})>r_u/u(p)-\delta/u(p)$ and by definition of $a$, $\delta/u(p)\leq\delta/a$. Therefore

 4) $d(\bar{c},\bar{m})>r_u/u(p)-\delta/a$.

Consider the case $k>0$. Since $d(\bar{c},\bar{x_l}) < 2d_k<D_k/2 $ and $d(\bar{c},\bar{x_l'}) < 2d_k<D_k/2 $, it follows  that $d(\bar{c},\bar{m})<2d_k$  . So by 4) $2d_k > r_u/u(p)-\delta/a$, from which it follows by 1) that $r_u/u(p)<2d_k + s/4 $. Therefore the supposition of statement 2) is satisfied. So for $k>0$, geodesics in $T$ have length less than $\epsilon/2$. For the case $k\leq0$, we have the apriori condition that geodesics in $T$ have length less than $\epsilon/2$.

From 3) $u(p)d(x_l,p)<r_u+\delta$ and $u(p)d(x_l',p)<r_u+\delta$. So 

$\displaystyle{d(x_l,p)<r_u/u(p)+\delta/u(p)< r_u/u(p)+\delta/a}$. Therefore 
$d(\bar{x_l},\bar{c}) < r_u / u(p) +\delta/a$ and similarly $d(\bar{x_l'},\bar{c}) < r_u/u(p)+\delta/a$. 

As $\bar{m}$ is the midpoint of $\bar{x_l}\bar{x_l'}$, it follows that 
$d(\bar{c},\bar{m})<d(\bar{c},\bar{x_l})$ or $d(\bar{c},\bar{m})<d(\bar{c},\bar{x_l'}).$ So $d(\bar{c},\bar{m})<r_u/u(p)+\delta/a$. 

This together with 4) implies that $\bar{m} \in T $. Without loss of generality, we may assume that $d(\bar{c},\bar{m}) \leq d(\bar{c},\bar{x_l})$. Thus $\bar{x_l} \in T$. So the geodesic segment $\bar{m}\bar{x_l}$ lies in $T$ and, consequently, has length less than $\epsilon/2$. So the length of the geodesic segment $\bar{x_l}\bar{x_l'}$ is less than $\epsilon$. It follows that if $l$ and $l'$ are greater than $M$, then $d(x_l,x_l') < \epsilon$. So $\{x_l\}$ is a Cauchy sequence. Since $X$ is complete the sequence converges to some $q_u$ in $X$.

If $\{y_l\}$ is another sequence  for which
$\limm{y_l}{P'}=r_u$ then by interlacing this sequence with $\{x_l\}$ we see that this new sequence also converges to $q_u$. So $q_u$ is the unique point in $X$ for which $r_u=\distt{p'}{q_u}{p'}{p' \in P'}$.
This completes the proof.

\end{proof}

For arbitrary $u$, $r_u$ is called the baryradius of $P$ relative to $u$, and $q_u$ is called the barycenter of $P$ relative to $u$.
When $u\equiv 1$, the baryradius of $P$ relative to $u$ is called the circumradius  of $P$ and is denoted by $r_1$, and the barycenter of $P$ relative to $u$ is called the center of $P$ and is denoted by $q_1$. Theorem \ref{T:center} for this case can be reformulated as follows: there exists a unique smallest closed ball containing $P$, the radius and center being $r_1$ and $q_1$, respectively. This result has been known for some time; see (\cite{bh}, Prop 2.7, P.179) for a proof. Furthermore, in this case, the $D_k/4$ condition can be replaced by $D_k/2$. To what extent this is so for other cases is not known to this author. The proof of Proposition 2.7  in \cite{bh} motivated the proof of Theorem \ref{T:center}. 

Let $(X,d)$ be a $\cat$ space and $P$ be a subset of $X$ of diameter less than $D_k/4$. Let $u:P\rightarrow \mathbb R$ be a non-negative bounded real valued function. Choose  $\lambda>0 $. Let $u':P \rightarrow \mathbb R$  be the function defined by setting $u'\equiv\lambda u$. Then   

\begin{eqnarray*}
r_{u'} &=& \inff{x}{\lambda u(p)d(x,p)}{p}{P}\\
    &=&  \lambda \inff{x}{ u(p)d(x,p)}{p}{P}\\
    &=& \lambda \distt{p}{q}{p}{p \in P}=\lambda r_u,
    \end{eqnarray*} 
We have established the following corollary.

\begin{corollary}\label{C:homo} Let $(X,d)$ be a $\cat$ space and $P$ a bounded subset of $X$ of diameter less than $D_k/4$. Let $u:P\rightarrow \mathbb R$ be a non-negative bounded real valued function. Let $\lambda>0 $ and set $u'\equiv \lambda u$.  Then $r_{u'}=\lambda r_u$ and the barycenter of $P$ relative to $u'$ is the same as the barycenter of $P$ relative to $u$.
\end{corollary}

Later, we make use of the following obvious fact.

\begin{corollary}\label{C:PtoP'} Let $(X,d)$ be a $\cat$ space and $P$ a bounded subset of $X$ of diameter less than $D_k/4$. Let $u:P\rightarrow \mathbb R$ be a non-negative bounded real valued function. If $P'$ is a subset of $P$ and $u\equiv 0$ on $P- P'$, then the barycenter of $P$ relative to $u$ is the same as the barycenter of $P'$ relative to $u$.
\end {corollary}

\centerline {
\includegraphics[width=2.5in]{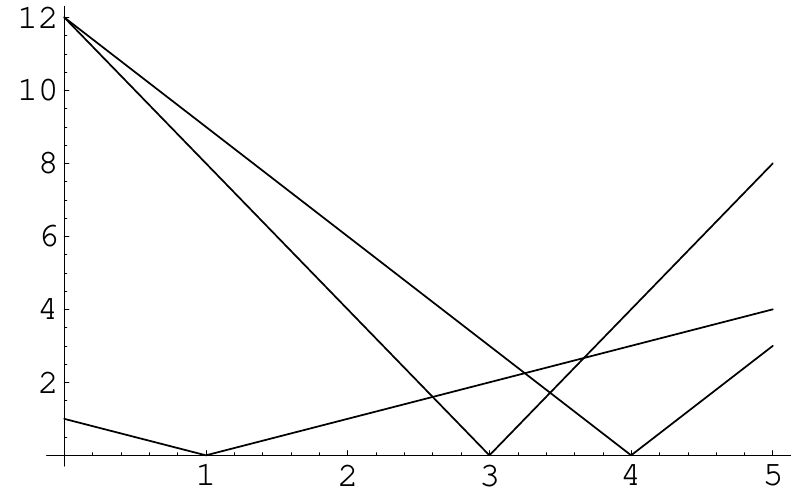}}

\begin{center}
Figure 2.
\end{center}

\begin{example} Let $P=\{1,3,4 \}\subset\mathbb{R}$. Set $u(1)=1$, $u(3)=4$ and $u(4)=3$. Then the barycenter of $P$ relative to $u$ is $x=13/4$, the $x$-coordinate of the lowest point on the graph  of $\max\{ |x-1|,4|x-3|,3|x-4| \mid x \in \mathbb{R}\}$.
\end{example}

\begin{example} Let $P=[0,1]\subset \mathbb{R}$ and $u:P\rightarrow\mathbb{R}$
be defined by $u(x)=x^2$. Then an easy calculation shows that the barycenter of $P$ relative to $u$ is $x=0.894$ to three decimal places.

\end{example} 
\begin{proposition}\label{P:close} Let $(X,d)$ be a $\cat$ space, $P $ a bounded subset of $X$ of diameter less than $D_k/4$. Let $u:P\rightarrow \mathbb R$ be a non-negative bounded real valued function. Suppose that $x$ is in $X$ and $P\subset \oball{x}{r}$ where $r$ is less than $D_k/2$. Then the barycenter $q_u $ of $P$ relative to $u$ lies in $\overline{B(x,r)}$. In particular, if $r=r_1$ is the circumradius of $P$ and $x=q_1$ is the circumcenter of $P$, then $q_u$ is in $\overline{B(x,r)}$.

\end{proposition}

\begin{proof} Suppose that $q_u$ is not in $\overline{B(x,r)}$. Set $d(x,q_u)=s$. Let $p$ be in $P$ and consider the geodesic triangle $\tri{p}{x}{q_u}$. It follows from Lemma \ref{L:close} and the triangle inequality that its perimeter is less than $2D_k$. So there is a  comparison triangle $\tribar{p}{x}{q_u}$ in $\model$. Let $y$ be the point on the geodesic segment $xq_u$ that satisfies $d(x,y)=r$, and $\bar{y}$ the corresponding comparsion point on $\bar{x}\bar{q_u}$. Set $d(x,q_u)=s$. By Lemma \ref{L:inequal2} there exist a positive number $\m{r}{s}$ such that 

\begin{equation*}
d(\bar{p},\bar{q_u})- d(\bar{p},\bar{y}) \geq \m{r}{s}.
\end{equation*} So \begin{equation*}
d(\bar{p},\bar{y})\leq d(\bar{p},\bar{q_u})-\m{r}{s}.
\end{equation*}
By the $\cat$ inequality, 

\begin{equation*}
d(p,y)\leq d(\bar{p},\bar{q_u})-\m{r}{s}=d(p,q_u)-m_k(r,s).
\end{equation*}
Therefore, for each $p$ in $P$

\begin{equation*}
u(p)d(p,y)\leq u(p)d(p,q_u)-u(p)\m{r}{s}.
\end{equation*}
By Lemma \ref{L:P'} we may assume that there exists a subset $P'$ of $P$, and  an $a>0$ such that $a=\inf\{u(p) \mid  p\in P' \}$ and, for $p$ in $P-P'$, $u(p)d(y,p)<r_u$, and
\begin{equation} \distt{p}{y}{p}{p \in P'}\geq r_u.\end{equation} Then for each $p$ in $P'$

\begin{equation*}
u(p)d(p,y)\leq u(p)d(p,q_u)-u(p)\m{r}{s} \leq u(p)d(p,q_u)-a\m{r}{s}.
\end{equation*} 
So
\begin{equation*}
\distt{p}{y}{p}{p \in P'} \leq r_u -a\m{r}{s}<r_u,\end{equation*} contradicting (1). So $q_u$ is in $\overline{B(x,r)}$. The second statement of the proposition follows immediately.

\end{proof}

Let $q_1$ be the circumcenter of a bounded set $P$ in $\mathbb R^n$. Let $q_u$ be the barycenter of $P$ relative to $u:P\rightarrow \mathbb{R}$. It follows from the proposition above and Jung's Theorem that  $\displaystyle{d(q_u,q_1) \leqq \sqrt{\frac{n}{2(n+1)}}diameter(P)}$. Analogous results can be established for general $\cat$ spaces by using the above proposition and estimates due to Lang and Schroeder \cite{ls}.

\begin{proposition}\label{P:convex} Let $P$ be a bounded subset in $C$, a closed convex subset of a complete $\catt$ space $X$ and let $u: P \rightarrow {\mathbb{R}}$ be bounded non-negative function from $P$ to $\mathbb{R}$. Then the barycenter of $P$ relative to $u$ lies in $C$.

\begin{proof} Let $x$ be in $X-C$, $p$ be in $P$ and $\pi(x)$ be the point in P such that $d(x,P)=d(x,\pi(x)$. Then by (\cite{bh}, Proposition 2.4, P.176) the Alexandroff angle $\angle_{\pi{x}}(x,p)\ge\pi/2$. The result follows by applying (\cite{bh}, Proposition 1.7, 
P.161) to an appropriate comparison triangle.

\end{proof}

\end{proposition}

It follows from the Proposition above, by setting $P=C$, and Example \ref{E:space}, Part 2,  that the barycenter is an intrinsic property of bounded $\catt $ spaces.

\begin{lemma}\label{L:contru} Let $P$ and $X$ be two non-empty sets and
let $w:X \times P\rightarrow \mathbb R$ and $w':X \times P\rightarrow \mathbb R$ be two non-negative bounded real valued functions. Set $r_w=\inff{x}{w(x,p)}{p}{P}$ and
$r_{w'}=\inff{x}{w'(x,p)}{p}{P}$. If for each $x$ in $X$ and each $p$ in $P$, $\vert w(x,p)-w'(x,p) \vert<\delta$, then $\vert r_w-r_{w'} \vert \leq \delta$.

\end{lemma}

\begin{proof} Since both $w$ and $w'$ are bounded functions $r_w$ and $r_{w'}$ are finite. Pick $\delta_1>0$. As $r_w=\inff{x}{w(x,p)}{p}{P}$, there exists
an $x$ in $X$ such that $\sup\{ w(x,p) \mid p \in P  \}< r_w+\delta_1$. For $p$  in $P$, $w'(x,p)<w(x,p)+\delta$. So

\begin{equation*}
r_{w'} \leq \sup\{ w'(x,p) \mid p \in P  \} \leq  \sup\{ w(x,p) \mid p \in P\} +\delta  \leq r_w +\delta + \delta_1
\end{equation*}
This being true for each $\delta_1 > 0$, it follows that $r_{w'} \leq r_w + \delta$. By reversing the role of $w$ and $w'$ above, we conclude that $r_{w}\leq r_{w'}+\delta$. So $\vert r_w-r_{w'}\vert \leq \delta$.

\end{proof}

\begin{definition}\label{D:fclose}
Let $(X,d)$ be a metric space, let $P$ be a non-empty subset of $X$ and let $\delta > 0$. Then a $\delta$-function, denoted by $\fdelta$, is a function from $P$ to $X$ such that $d(p,\fdelta(p))<\delta$, for every $p$ in $P$.
\end{definition}

Next we show that if the points $P'$  are sufficiently close to $P$ and the function $u'$ is sufficiently close to $u$, then the barycenters and baryradii of the respective sets and functions will be close.

\begin{theorem}\label{T:con} Let $X$ be a complete $\cat$ space and let $P$ be a subset of $X$ with diameter less than $D_k/4$. Let $u:P \rightarrow \mathbb R$ and $u':P\rightarrow \mathbb R$ be non-negative and bounded real valued functions. Let $q_u$ be the barycenter of $P$ relative to $u$, and let $r_u$ be the baryradius relative to $u$. Let $q_{u'}$ be the barycenter of $P$ relative to $u'$, and let $r_{u'}$ be the baryradius relative to $u'$.

\begin{enumerate}
\item 
For every $\epsilon >0$, there exists a $\delta>0$ such that if $\vert u(p)-u'(p)\vert<\delta$, for every $p$ in $P$, then $\vert r_u-r_{u'}\vert<\epsilon$.
\item Given $u$ and $\epsilon >0$, there exists a $\delta>0$ such that if $\vert u(p)-u'(p)\vert<\delta$, for each $p$ in $P$, then $d(q_u,q_{u'}) < \epsilon$.
\item Given $u$ and $\epsilon >0$, there exists a $\delta > 0$ so that if $\fdelta$ is a $\delta$-function from $P$ to $X$, then the barycenter of  $P$ relative to $u$ and the baryradius of $\fdelta(P)$ relative to $u$  \footnote{ Here we assume that $u(p')= \sup\{u(p) \mid f_\delta (p)=p' \}$. } are within $\epsilon$ of $q_u$ and $r_u$, respectively.
\end{enumerate}
  \end{theorem}
  
  \begin{proof} \setcounter{equation}{0} By Proposition \ref{P:close}, both $q_u$ and $q_{u'}$  are in $\overline{B(q_1,r_1)}$ where $r_1$ is the circumradius of $P$ and $q_1$ is its circumcenter. So
  
  \begin{equation}  r_u=\inff{x}{u(p)d(x,p)}{p}{P}=\inff{x \in \overline{B(q_1,r_1)}}{u(p)d(x,p)}{p}{P}
  \end{equation}
  and 

 \begin{equation}
  r_{u'}=\inff{x}{u'(p)d(x,p)}{p}{P}=\inff{x \in \overline{B(q_1,r_1)}}{u'(p)d(x,p)}{p}{P}.
  \end{equation}

\begin{paragraph}{Part 1.}  Let $\epsilon>\epsilon'>0$. Set $\delta=\epsilon'/2r_1$. If  $\vert u(p)-u'(p) \vert <\delta$ for every $p$ in $P$, then for every $p$ in $P$ and $x \in \oball{q_1}{r_1},$   
   \begin{equation*}
  \vert u(p)d(x,p)-u'(p)d(x,p) \vert = \vert u(p)-u'(p)\vert d(x,p)<\frac{\epsilon'}{2r}2r = \epsilon'.
  \end{equation*}
  By (1),(2) and Lemma \ref{L:contru}, $\vert r_u-r_{u'}\vert \leq \epsilon'<\epsilon$. This completes the proof of Part 1.
 \end{paragraph}

\begin{paragraph}{Part 2.} Choose $\epsilon > 0$. Set $D=\oball{q_1}{r_1} - \ball{q_u}{\epsilon} $ and
\begin{equation}
 s=\inff{x \in D}{u(p)d(x,p)}{p}{P}.
 \end{equation}
 By the uniqueness of $q_u$, $ s-r_u >0$. Set
\begin{equation} 
\delta<(s-r_u)/4r_1\end{equation} and suppose that $\vert u(p)-u'(p) \vert < \delta$ for each $p$ in $P$. Then for each $x$ in $\oball{q_1}{r_1}$ and each $p$ in $P$,
 \begin{equation*}
 \vert u(p)d(x,p)-u'(p)d(x,p)\vert =\vert u(p)-u'(p)\vert d(x,p) < 2\delta r_1.
 \end{equation*}  By Lemma \ref{L:contru} 
 \begin{equation*}
 \vert \inff{x \in D}{u(p)d(x,p)}{p}{P}-\inff{x \in D }{u'(p)d(x,p)}{p}{P} \vert \leq 2 \delta r_1;
\end{equation*} so by (3)

\begin{equation}
\vert s-\inff{x \in D}{u'(p)d(x,p)}{p}{P} \vert \leq 2 \delta r_1.
\end{equation} Also by Lemma \ref{L:contru},
\begin{equation*}
\vert \inff{x \in \ball{q_u}{\epsilon} }{u(p)d(x,p)}{p}{P}-\inff{x \in \ball{q_u}{\epsilon} }{u'(p)d(x,p)}{p}{P} \vert \leq 2 \delta r_1,
\end{equation*} or 
 \begin{equation}
 \vert r_u-\inff{x \in \ball{q_u}{\epsilon}}{u'(p)d(x,p)}{p}{P} \vert < 2 \delta r_1.
 \end{equation}
 By 5)
 \begin{equation*}
 s-2\delta r_1 \leq \inff{x \in D }{u'(p)d(x,p)}{p}{P},
 \end{equation*}
 and by 6)
 \begin{equation*}
 \inff{x \in \ball{q_u}{\epsilon} }{u'(p)d(x,p)}{p}{P} \leq r_u+2\delta r_1. 
 \end{equation*}
 By 4)
 
  $r_u+2\delta r_1 < s-2\delta r_1$, so 
 \begin{equation*}
 \inff{x \in \ball{q_u}{\epsilon} }{u'(p)d(x,p)}{p}{P} < \inff{x \in D}{u'(p)d(x,p)}{p}{P}.
 \end{equation*}
Therefore $q_{u'}$, the barycenter of $P$ relative to $u'$, is in $\ball{q_u}{\epsilon} $. This completes the proof of Part 2.
 
 \end{paragraph}

\begin{paragraph}{Part 3.} Let $\epsilon>0$ and set $D=\oball{q_1}{r_1} - \ball{q_u}{\epsilon} $. By Corollary \ref{C:homo}, we may assume that $\vert u(p) \vert \leq 1$ for each $p$ in $P$. Set $s=\inff{x \in D}{u(p)d(x,p)}{p}{P}$, and choose $\delta$ so that 
\begin{equation}
s-r_u>2\delta. 
\end{equation}
 Suppose that $f_{\delta}$ is a $\delta$-function. Then for each $x$ in $X$ and each $p$ in $P$, 
\begin{equation*}
\vert u(p)d(x,p)-u(p)d(x,\fdelta(p)) \vert < \delta. 
\end{equation*}
Then by Lemma \ref{L:contru}

\begin{equation*}
\vert \inff{x \in D}{u(p)d(x,p)}{p}{P}-\inff{x \in D}{u(p)d(x,\fdelta(p))}{p}{P} \vert \leq  \delta, 
\end{equation*}
So\begin{equation}
\vert s-\inff{x \in D}{u(p)d(x,\fdelta(p))}{p}{P} \vert \leq  \delta. 
\end{equation}
Also by Lemma \ref{L:contru} 

\begin{equation*}
\vert \inff{x \in \ball{q_u}{\epsilon}}{u(p)d(x,p)}{p}{P}-\inff{x \in \ball{q_u}{\epsilon}} {u(p)d(x,\fdelta(p))}{p}{P} \vert \leq  \delta, 
\end{equation*}
or

\begin{equation}
\vert r_u-\inff{x \in \ball{q_u}{\epsilon}} {u(p)d(x,\fdelta(p))}{p}{P} \vert \leq \delta. 
\end{equation}
Therefore,

\begin{equation*}
s-\delta < \inff{x \in D}{u(p)d(x,\fdelta(p))}{p}{P}
\end{equation*} 
and 

\begin{equation*}
\inff{x \in \ball{q_u}{\epsilon}} {u(p)d(x,\fdelta(p))}{p}{P} \leq r_u +\delta.
\end{equation*}
 By 7) $r_u +\delta < s-\delta$. Therefore the barycenter of $\fdelta$, relative to $\fdelta(P)$, is in $\ball{q_u}{\epsilon}$ and this completes the proof of the first half of Part 3). The argument for the last half is essentially the same as the proof of Part 1) and will be omitted.
\end{paragraph}
\end{proof}

Suppose that $\{ X_n \}$ is a uniformly bounded sequence of complete $\catt$ spaces that converge to the complete
$\catt$ space $X$. We will define what it means for a sequence of uniformly bounded non-negative functions $\fn{u_n}{X_n}{\bld}$ to converge to a non-negative function $\fn{u}{X}{\bld}$. Then we will show that $\underset{n \to \infty}{\lim}{q_n=q_u}$ where $q_n$ is the barycenter of $X_n$ relative to $u_n$ and $q_u$ is the barycenter of $X$ relative to $u$.

Let $X$ be a subset of a metric space $Z$. Then $U_{\epsilon}(X)=\{ z \in Z \vert d(z,X)< \epsilon \}$ is called the 
$\epsilon$-neighborhood of $X$ in $Z$.

Let $X$ and $Y$ be subspaces of a metric space $Z$. Then $d_H(X,Y)$, the Hausdorff distance between $X$ and $Y$, is defined by $d_H(X,Y)= \inf\{\epsilon \vert X\subseteq U_{\epsilon}(Y) $ and $  Y\subseteq U_{\epsilon}(X) \}$. A
$\delta$-resolution $\R(X,Y,Z,\delta)$ of $X$ and $Y$, or resolution, consists of a metric space $Z$ and isometries $\fn{i_1}{X}{Z}$ and 
$\fn{i_2}{Y}{Z}$ such that $d_H(i_1(X),i_2(Y))<\delta$. 

\begin{itemize}
\item For $x$ in $X$ and $y$ in $Y$ we will identify $x$ with $i_1(x)$ and $y$ with $i_2(y)$. We identify $X$ with $i_1(X)$ and $Y$ with $i_2(Y)$.
 \end{itemize}

Let $X$ and $Y$ be metric spaces. The Gromov-Hausdorff distance between $X$ and $Y$ is given by
$d_{GH}(X,Y)=\inf \{ \delta \vert \R(X,Y,Z,\delta) \text{ is a $\delta$-resolution of $X$ and $Y$} \}$.

A sequence of metric spaces $ \{X_n \}$ converges to a metric space $X$ provided that 
$\underset{n\rightarrow \infty}\lim d_{GH}(X_n,X)=0$. For more  information about the Gromov-Hausdorff metric,
see (\cite{bh},\cite{bbi}).

For each natural number $n$, let $\fn{u_n}{X_n}{\bld}$ be a real valued function defined on metric space $X_n$ and let  $\fn{u}{X}{\bld}$ be 
a real valued function defined on metric space $X$. We say that $\{u_n \}$ converges to $u$ provided the following
two conditions are satisfied. There exists a sequence of resolutions $\R(X_n,X,Z_n,\delta_n)$ such that 
$\underset{n\rightarrow \infty}\lim \delta_n=0$, and for every $\epsilon>0$, there exist an $N>0$ such that if $n>N$, $x_n \in X_n \subseteq Z_n$, $x \in X  \subseteq Z_n$ and $d(x_n,x)<\delta_n$, then $d(u_n(x_n),u(x))<\epsilon$. Implicit in this definition is the fact that $X_n$ converges to $X$.

For each natural number $n$, let $q_n$ be a point in the metric space $X_n$ and let $q$ be a point in the metric space 
$X$. Then we say that the sequence $\{ q_n \}$ converges to $q$ provided that the following two conditions are satisfied. There exists a sequence of resolutions $\R(X_n,X,Z_n,\delta_n)$ such that 
$\underset{n\rightarrow \infty}\lim \delta_n=0$, and for every $s>0$, there exists an $N>0$ such that if $n>N$ then the
ball $\ball{q}{s}$ in $Z_n$ contains $q_n$.

\begin{lemma}\label{L:fap} Let $X$ and $Y$ be bounded metric spaces such that $b=\max\{ dia(X),dia(Y)\}$. Let 
$\fn{f}{X}{\bld}$ and $\fn{g}{Y}{\bld}$ be non-negative real valued functions such that 
$m=\max\{ \underset{x}\sup{f(x)},\underset{y}\sup{g(y)}  \}$. Let $ \R(X,Y,Z,\delta) $ be a $ \delta$-resolution of $X$ and $Y$ and let $\epsilon$ be a positive number. Assume that if $x$ is in $X$, $y$ is in $Y$ and $d(x,y) < \delta$,
then $\vert f(x)-g(y)  \vert<\epsilon$. 

If $x$ and $x'$ are in $X$, $y$ and $y'$ are in $Y$ so that $d(x,y)<\delta$ and  $d(x',y')<\delta$, then 
$g(y')d(y,y')-\delta' \leq f(x')d(x,x') \leq g(y')d(y,y') + \delta'$, where $\delta'\leq 2\delta m +\epsilon b +2\delta \epsilon$.
\end{lemma}

\begin{proof} Several applications of the triangle inequality imply that $d(y,y')-2\delta < d(x,x') <d(y,y') + 2\delta$. Since 
$d(x,x')>0$ there exists $\delta''=\delta''(y,y') \leq 2\delta$ such that

1) $0 <d(y,y')-\delta''< d(x,x') < d(y,y') + 2\delta$. 

Since $\vert f(x')- g(y') \vert < \epsilon$, 

2) $g(y')-\epsilon <f(x') < g(y') + \epsilon$.

It follows from 1) and 2) that

$(g(y')-\epsilon)(d(y,y')-\delta'')< f(x')d(x,x') < (g(y') +\epsilon)(d(y,y') + 2\delta)$

So

$ g(y')d(y,y') - (\delta'' g(y')+\epsilon d(y,y') - \epsilon \delta) \leq  f(x')d(x,x') \leq  g(y')d(y,y') + (2 \delta g(y') + \epsilon 
d(y,y') +2 \delta \epsilon) $. 

To complete the proof, set $\delta' = 2 \delta g(y') + \epsilon d(y,y') +2 \delta \epsilon $ which is equal to or less than $2 \delta m + \epsilon b + 2\delta \epsilon$.

\end{proof}

We define a sequence of metric spaces $\{ X_n \}$ to be uniformly bounded if there exists a positive number $b$ such that
the $dia(X_n)<b$ for each $n$, and we define a sequence of real valued functions $\fn{u_n}{X_n}{\bld}$ to be uniformly bounded if there exists a positive number $m$ such that $\sup\{\vert u_n(x)\vert\, \vert x\in X_n\, \text{and}\,  n\in \mathbb N \} \leq m$.

\begin{theorem}\label{T:aprox}  Suppose that $X_n$ is a sequence of uniformly bounded complete $\catt$ spaces 
that converge to the complete $\catt$ space $X$, and $\fn{u_n}{X_n}{\bld}$  is a sequence of uniformly bounded non-negative functions that converge to the non-negative function $\fn{u}{X}{\bld}$. Then  $\underset{n \to \infty}{\lim}{q_n=q_u}$ where $q_n$ is the barycenter of $X_n$ relative to $u_n$ and $q_u$ is the barycenter of $X$ relative to $u$.

\end{theorem}

\begin{proof} By the remark following Proposition \ref{P:convex} the barycenter of each $X_n$ and $X$ is independent of the ambient space. Consider $\ball{q_u}{2s}$, the open ball in $X$ of radius $2s$ centered at $q_u$. Let $b$ be a uniform bound
on the spaces $X_n$ and let $m $ be a uniform bound on the functions $u_n$. Let

$r=\inff{x \in X-\ball{q_u}{s}}{u(x')d(x,x')}{x'}{X}-r_u$, where $ r_u=\inff{x \in X}{u(x')d(x,x')}{x'}{X}$, the baryradius of $X$ relative to $u$.

The barycenter being unique implies that $r>0$.  Since $u_n$ converges to $u$, there exists a sequence of resolutions
$\R(X_n,X,Z_n,\delta_n)$ such that $\underset{n\rightarrow \infty}\lim \delta_n=0$. Choose a positive $\epsilon$ so that $\epsilon< \min{\{r/4,r/{4b}\}}$. Let $t<s/2$ be chosen so that if $x$ is in $B(q_u,t)$ then 
$\sup{u(p)d(x,p) \mid <r_u+r/4}$. Furthermore, we can pick an $N>0$ such that 

1) if $x_n \in X_n \subset Z_n$, $x  \in X \subset Z_n$ and $d(x_n,x)<\delta_n$, then $\vert u_n(x_n)-u(x) \vert< \epsilon$,

2) $ \delta_n<t$, and

3) $2\delta m+\epsilon b +2\epsilon \delta < r/4$.

By 1), 3) and Lemma \ref{L:fap} we may assume that if $x_n$ and $x_n'$ are in $X_n$, $x$ and $x'$ are in $X$, $d(x_n,x)<\delta_n$ and $d(x_n',x')<\delta_n$, then

4) $u(x')d(x,x')-r/4< u_n(x_n')d(x_n,x_n') <u(x')d(x,x')+r/4$.

From the definition of a resolution, we have

5) $d_H(X_n,X)<\delta_n$.

So there exists an $x_n$ in $X_n$ such that $d(x_n,q_u)<\delta_n<t$. Pick an arbitrary $x_n'$ in $X_n$. Again by 5) there 
exists an $x'$ in $X$ such that $d(x_n',x')<\delta_n$. Therefore by 4) $u_n(x_n')d(x_n,x_n')<u(x')d(q_u,x')+r/4 \leq r_u +r/4$. 
Therefore, by the inequality in 4) with $n>N$, we have

6) $x_n $ is in $\ball{q_u}{s}$ and $\underset{x_n' \in X_n}\sup \{u(x_n')d(x_n,x_n')\} \leq r_u+r/4$.

Now suppose that $z_n \in X_n-\ball{q}{2s}  $. By 5) there exists a $z$ in $X$ such that $d(z_n,z)<\delta_n<t$. On account of the definition of $t$, we conclude that $z$ is not in $B(q_u,s)$.
So there exists a $z'$ in $X$ such that $u(z')d(z,z')>r_u +r$. By 5), once again, there is an
$x_n'$ in $X_n$ such that $d(x_n',x')<\delta_n$.

By the first inequality in 4) and the inequality above $r_u +r-r/4<u(z')d(z,z')-r/4<u_n(z_n')d(z_n,z_n')$. Therefore,

$\sup\{ u_n(z_n')d(z_n,z_n' \}\geq r_u +3/4r$.
 
Since $z_n$ was an arbitrary point in $X_n-\ball{q_u}{2s}$, it follows from 6) and the inequality above that $q_n $ is in $\ball{q_u}{2s}$. Therefore, $\underset{n \rightarrow \infty}\lim q_n=q_u$.

\end{proof}

More can be said about the sequence $\{ q_n\}$. Let $r>0$. It is easy to see that $d_{GH}(\overline{B(q_n,r)},\overline{B(q_u,r)}) \leq d_H(\overline{B(q_n,r)},\overline{B(q_u,r)})< d(q_n,q_u) + 2\delta_n$. Since $\underset{n \rightarrow \infty}\lim{d(q_n,q_u) + 2\delta_n}=0$, it follows that the pointed sequence $ (X_n,q_n)$ converges to $(X,q_u)$ (See P.76,\cite{bh}).

As in Definition  \ref{D:disttt}, let $X$ be a $\cat$ space, and $P$ be a subspace of  $X$ with diameter that is finite and less than $D_k/4$. Let $u:P\rightarrow \mathbb R$ be a bounded and non-negative real valued function from $P$ to $\mathbb R$. Let $t>0$ and for each  $x$ in $X$ define $r(x)=\sup \{u(p)d(x,p)^t \mid p \in P \}$. Set $r_u^t=\dis{x}{u(p)d(x,p)^t}$. Then $r(x)^{1/t}=\sup\{ u(p)^{1/t}d(x,p) \mid p \in P\}$ and $r_{u^{1/t}}= \underset{x}\inf\{ r(x)^{1/t}\}=\dis{x}{u^{1/t}(p)d(x,p)}$. By Theorem \ref{T:center} there exists a unique $q$ in $X$ such that $r_{u^{1/t}}= \sup\{ u(p)^{1/t}d(q,p) \mid p \in P \}$. It follows that
$r_u^t=(r_{u^{1/t}})^t=\sup\{u(p)d(q,p)^t \mid p \in P\}$. Thus we see that the distance function $d(x,y)$ in Theorem \ref{T:center} can be replaced by $d(x,y)^t$ where $ t $ is any positive number. Furthermore, the properties of the barycenter carry over to this more general environment.

A  barycenter, also called the center of mass, is defined in a slightly different manner in (\cite{bbi}, P. 334),
and this definition is equivalent to the definition of center of mass, with Dirac measure, in \cite{Jost}.When compared to the definition in this paper, sum replaces supremum. Their definition provides a more sensitive, but less general concept.

\section{Applications}\label{S:retract} 

Let $ f:X \rightarrow X$ be a function from $X$ to itself. An element $x$ in $X$ is a fixed point of $f$ provided that
$f(x)=x$. The circumcenter has been used to establish fixed points of isometries in $\catt$ spaces (\cite{bh}, P.179). There is also a connection between fixed points and barycenters.

\begin{proposition}\label{P:fix1} Let $P$ be a bounded subset of a $\catt$ space $X$ and let $ f:X \rightarrow X$ be an isometry that is invariant on $P$. If $u:P\rightarrow \mathbb{R}$ is a bounded function such that $f(u(p))=u(p)$ for each $p$ in $P$, then the barycenter of $P$ relative to $u$ is a fixed point of $f$.

\end{proposition}

\begin{proof}
Let $q_u$ be the barycenter of $P$ relative to $u$. The proof follows immediately from the fact that
$\displaystyle{\distt{p}{q_u}{p}{p \in P}=\distt{f(p)}{f(q_u)}{f(p)}{p \in P}}$ and the uniqueness of the barycenter.

\end{proof}

A similar proof establishes the next proposition.

\begin{proposition}

 Let $X$ be a $\catt$ space and let $\Gamma$ be a group of isometries acting on $X$. Let $P$ be a bounded set of orbits and let $u:P \rightarrow \mathbb{R}$ be a bounded function such that $u(p)=u(p')$
 when $p$ and $p'$ are in the same orbit. Then the barycenter of $P$ relative to $u$ is a common fixed point of $\Gamma$.
 
\end{proposition} 

And similar results can be formulated in $\cat$ spaces.

Let $(Y,d)$ be a metric space and $X$ be subset of $Y$. Then $X$ is a retract of $Y$ if there exists a retraction $f$ of $Y$ onto $X$. That is a continuous function $f:Y \rightarrow X$ such that $f(x)=x$ for every $x$ in $X$. 

Let $(X,d)$ be a metric space. Then $X$ is an absolute neighborhood retract (ANR) if and only if for every embedding, that is a homeomorphism, $h$ of $X$ into a metric space $Y$ as a closed subspace, there exists an open neighborhood $U$ of $h(X)$ in $Y$ and a retraction $f$ of $U$ onto $h(X)$. If it is always possible to choose the retraction from $Y$ to $h(X)$, then $X$ is called an absolute retract (AR). The reader is referred to \cite{bor,gd,hu} for a detailed development of this subject. We will make use of the following two theorems.

\begin{definition}\label{d:k}
Let $(X,d)$ be a geodesic metric space. If for each $x$ in $X$ there exists an 
$r_x>0$ such that the the ball $B(x,r_x)$, with the induced metric, is a $\cat$
space, then $X$ is said to be of curvature $\leq k$.
\end{definition}

It follows from a theorem of  Hu (\cite{hu},IV,Thm 4.1) that  a complete space of curvature $ \leq k $ is an absolute neighborhood retract (\cite{bh} P.209-210;\cite{fg}, Remark 4.2).
We give a different proof using barycenters.

\begin{theorem}(Arens-Eells)\label{T:arens} Every metric space can be isometrically embedded as a closed subspace of some normed linear space.

\end{theorem}
See (\cite{gd}, P.597) for a proof of this result. For a proof of the following theorem see (\cite{gd}, P.164, P.280).

\begin{theorem}\label{T:aranr} \item a) A metric space $(X,d)$ is an AR if and only if there exists an embedding $h$ of $X$ into a normed linear space $N$ and a  retraction $r:N\rightarrow h(X)$ of $N$ onto $h(X)$.
\item b) A metric space $(X,d)$ is an ANR if and only if there exists an embedding $h$ of $X$ into a normed linear space $N$, a neighborhood $U$ of $h(X)$ and a retraction $r:U\rightarrow h(X)$ of $U$ onto $h(X)$.
\end{theorem}
 The next lemma can be found in the body of the proof of the Dugundji extension theorem (\cite{gd}, P.163-164).
\begin{lemma}\label{L:howfar} Let $X$ be a closed subspace of a metric space $(Y,d)$. Suppose that $x$ and  $z$ are in $X$, and $y$ and $v$ are in $Y-X$.  Suppose that $y$ is an element of $\overline{B(v,\frac{d(v,X)}{2})}$  and $d(v,z)<2d(v,X)$. Then $d(x,z)< 6d(x,y)$. 
\end{lemma}

\begin{proof} Several applications of the triangle inequality yield
\begin{eqnarray*}
d(x,z) &\leq& d(x,y)+d(y,v)+d(v,z)\\
       &<& d(x,y)+\frac{d(v,X)}{2}+2d(v,X)\\
       &<& d(x,y)+\frac{5}{2} d(v.X)  
       \end{eqnarray*}
Now 
\begin{equation*} d(v,X) \leq d(v,y)+d(y,X) \leq d(y,x)+d(y,x)=2d(y,x).
\end{equation*} 
Therefore, $d(x,z)< 6d(y,x)$.      

\end{proof}

Let $(X,\script{T})$ be a topological space and let $\script{U}$ be an open cover of $X$. Then $\script{V}$ is an open refinement of $\script{U}$ provided that $\script{V}$ is a cover of $X$ consisting of open sets, and each $V$ in $\script{V}$ is a subset of some member $U$ of $\script{U}$.  Now $\script{V}$ is a neighborhood finite open refinement of $\script{U}$ provided that $\script{V}$ is an open refinement of $\script{U}$ and, for each $x$ in $X$, there exists a neighborhood $N(x)$ of $x$ that meets only finite many $V $'s in $\script{V}$. If $A$ is a subset of $X$ then the $\starr{(A,\script{V})}= \underset{V\in \script{V}}\bigcup\{V \mid V\cap A \ne \emptyset \}$. Then $\script{V}$ is an open star refinement of $\script{U}$ provided that, each member of $\script{V}$ is an open subset of $X$, and $\{ \starr(V,\script{V}) \mid V\in \script{V} \}$ is a refinement of $\script{U}$. If $X$ is a metric space, then every open cover of $X$ has a neighborhood finite open refinement and an open star refinement as well. See (\cite{d}, P.167-168).

Although Theorem \ref{T:ar} that follows is a consequence of Theorem \ref{T:anr}, we give an independent proof which is essentially the same as Dugundji's proof of his extension theorem (\cite{gd}, P.163-164) with convex sums replaced by barycenters.

\begin{theorem}\label{T:ar}Let $(X,d)$ be a complete  $\cat$ space with $k \leq 0$. Then $X$ is an absolute retract.
\end{theorem}

\begin{proof} By Theorem \ref{T:arens} we may assume that $X$ is a closed subspace of a normed linear space $N$, and the distance function $d(x,y)$ on $X$ agrees with the norm metric on $X$.

Cover $N-X$ with the set of  open balls $\script{B}=\{B(y,d(y,X)/2)\mid y \in N-X \}$. Let $\mathscr{C}=\{ V_\lambda \mid \lambda \in \Gamma\}$, where $\Gamma $ is an index set, be a neighborhood finite open refinement. 

For each $V_\lambda$, choose a ball $B(v_\lambda,d(v_\lambda,X)/2)$ that contains $V_\lambda$, and also choose an $x_\lambda$ in $X$ such that the distance from $v_\lambda$ to $x_\lambda$ is less than $2d(v_\lambda,X)$.

Choose the partition of unity corresponding to $\script{C}$ defined by 
\begin{equation*}\displaystyle{u_\lambda(x)=\frac{d(x,X-V_\lambda)}{\underset{\beta \in \Gamma}\sum d(x,X-V_\beta)}}, \text{for each $\lambda$ in $\Gamma$.}
\end{equation*} 

We define the retraction $f:N \rightarrow X$ as follows. If $x$ is in $X$ then $f(x)=x$. Suppose that $y$ is not in $N-X$. Let $\sett{ V_{\lambda_1}}{V_{\lambda_m}}$ be those members of $\script{C}$ that contain $y$. Set $f(y)$ equal to the barycenter of $x_{\lambda_1}, \dots , x_{\lambda_n}$ relative to  $ \sett{u_{\lambda_1}(y)}{u_{\lambda_m}(y)}$ as given by Theorem \ref{T:center}.

Let $\epsilon>0$. Suppose $y$ is in $N-X$. Since $\mathscr{C}$ is a neighborhood finite open refinement of $\script{B}$ there is a neighborhood $N(y)$ of $y$ that has non-empty intersection with only a finite number of members of $\mathscr{C}$. Choose the neighborhood so that if $V$ is in $\script{C}$ and $V\cap N(y)\neq \varnothing$ then $y$ is in $\bar{V}$. Let $V_{\lambda_1},\dots, V_{\lambda_n}$ be those members of $\script{C}$ whose closures have non-empty intersection with $N(y)$. 

If $y'$ is in $N(y)$, then it follows from the definition of $f$ and Corollary \ref{C:PtoP'} that $f(y')$ is the barycenter of $\sett{x_{\lambda_1}}{x_{\lambda_n}}$ relative to $\sett{u_{\lambda_1}(y')}{u_{\lambda_n}(y')}$. By Theorem \ref{T:con}, Part 2, there exists a $\delta'>0$ such that if $\vert u_{\lambda_i}(y)-u_{\lambda_i}(y') \vert<\delta'$ for each $i$ in $[n]$\footnote{$[n]=\sett{1}{n}$}, then  the barycenter relative to $\sett{u_{\lambda_1}(y')}{u_{\lambda_n}(y')}$ is within $\epsilon$ of the barycenter relative to $\sett{u_{\lambda_1}(y)}{u_{\lambda_n}(y)}$.  Since the members of the partition of unity are continuous, there exists a $\delta>0$ such that if $d(y,y')<\delta$ then, for each $i$ in $[n]$, $\vert u_{\lambda_i}(y)-u_{\lambda_i}(y') \vert<\delta'$. That is for this choice of $\delta$, if $d(y,y')<\delta$, then $d(f(y),f(y')<\epsilon$. So $f$ is continuous at $y$.
 
 Now suppose that $x$ is in $\bd(X)$ the boundary of $X$. Let $y$ be in $B(x,\delta)-X$ where $\delta<\epsilon/6$. Suppose $V_{\lambda_1},\dots, V_{\lambda_m}$ are those members of $\script{C}$ that contain $y$. By Lemma \ref{L:howfar} each $x_{\lambda_i}$ is in $B(x,\epsilon)$. By Proposition \ref{P:close} the barycenter of $x_{\lambda_1}, \dots, x_{\lambda_m}$
 relative to $\sett{u_{\lambda_1}(y)}{u_{\lambda_m}(y)}$ is in $B(x,\epsilon)$. Therefore, $f(y)$ is in $B(x,\epsilon)$. Thus $f$ is a continuous retraction of $N$ onto $X$. By Theorem \ref{T:ar}, Part a), $X$ is an absolute retract.
 
\end{proof}

\begin{theorem}\label{T:anr} Let $(X,d)$ be a complete metric space of curvature $\leq k$. Then $(X,d)$ is an absolute neighborhood retract.

\end{theorem}

\begin{proof}  As in the proof of Theorem \ref{T:ar} we may assume that $X$ is a closed subspace of a normed linear space $N$, and the distance function $d(x,y)$ on $X$ agrees with the norm metric on $X$.

Let $d_k/4$ be a positive number less than $D_k/4$. About each point $x$ in $\bd(X)$, the boundary of $X$, select an open ball $B(x,s_x)$ such that $s_x< d_k/4$. From the definition of curvature $\leq k$ and Example \ref{E:space}, Part 3, we may assume that $B(x,s_x)\cap X$ is a $\cat$ space. Let $W=X \cup \underset{x\in bd(X)}\bigcup B(x,s_x/18) $. Then $W$ is an open subset of $N$ containing $X$.

Let $\script{B}=\{ B(x,s_x/18)-X \mid x\in bd(X)\}$. Let $\script{C}$ be an open star refinement of $\script{B}$. For each $y$ in $W-X$ choose a member $U_y$ of $\script{C}$ that contains $y$. Then choose a ball $B(y,r_y) \subset U_y$ such that $r_y<d(y,X)/2.$ Set $\script{D}=\{ B(y,r_y) \mid y \in W-X \}$.    Let $\script{E}$ be a neighborhood finite open refinement of $\script{D}$. For each $V$ in $\script{E}$ choose a $B(y,r_y)$ in $\script{D}$ such that $V\subset B(y,r_y)$. Then choose an $x_V$ in $X$ such that $d(y,x_V) < 2d(y,X)$. We define a partition of unity as follows. For each $V$ in $\script{E}$ define  \begin{equation*}\displaystyle{u_V(y)=\frac{d(y,N-V)}{\underset{V'\in \script{E}}\sum d(y,N-V')}}.
\end{equation*}

Next we define a retraction $f:W\rightarrow X$. Of course $f(x)=x$ for $x$ in $X$. Let $y$ be in $W-X$. Since  $\script{E}$ is a neighborhood open finite refinement of
$\script{D}$, there exists a neighborhood $N(y)$ of $y$ and finitely many members $V_{\lambda_1} \dots V_{\lambda_n}$ of $\script{E}$ such that

\begin{itemize}

\item if $V$ is in $\script{E}$, then $N(y)\cap V\ne\emptyset$ if and only if $V=V_{\lambda_i}$, for some $i\in [n]$,

 and
 \item $y$ is in $\overline {V_{\lambda_i}}$ for each $i$.
\end{itemize}

For each $i\in [n]$ let $B(y_{\lambda_i},r_{\lambda_i})$ be the ball in $\script{D}$ associated with $V_{\lambda_i}$ and let $U_{\lambda_i}$ be the member of $\script{C}$ chosen to contain $B(y_{\lambda_i},r_{\lambda_i})$. Then $V_{\lambda_i}\subset B(y_{\lambda_i},r_{\lambda_i})\subset U_{\lambda_i}$. Since $\script{E}$ is a cover of $W-X$ at least one $V_{\lambda_i}$, say $V_{\lambda_1}$, contains $y$. Therefore $U_{\lambda_1}$ contains $y$. Since $\script{C}$ is a star refinement of $\script{B}$, there exists a $B(x,s_x/18)$
such that $\underset{i\in[n]}\bigcup V_{\lambda_i}\subset \underset{i\in[n]}\bigcup U_{\lambda_i}\subset B(x,s_x/18)$. Let $x_{V_{\lambda_i}}$ be the point in $X$ associated with $V_{\lambda_i}$, $i\in [n]$. From the definition of the $x_{V_{\lambda_i}}$ and Lemma \ref{L:howfar}, $d(x,x_{V_{\lambda_i}})<s_x/3$. Therefore each $x_{V_{\lambda_i}}$ is in $B(x,s_x/3)$ and $\dia{ \sett{x_{V_{\lambda_1}}}{x_{V_{\lambda_n}}}  }<2s_x/3<d_k/4$. Since $B(x,s_x)$ is a $\cat$ space, we define $f(y)$ to be the barycenter of 
$\sett{x_{V_{\lambda_1}}}{x_{V_{\lambda_n}}}$ relative to $\sett{u_{V_{\lambda_1}}(y)}{u_{V_{\lambda_n}}(y)}$. Theorem \ref{T:center} guarantees that $f(y)$ exists and is unique. An application of the triangle inequality shows that, for $i \in [n]$, if $x'$ is in $\overline{B(x_{V_{\lambda_i}},2s_x/3)}$, then $x'$ is in $B(x,s_x)$. Therefore, for each $i$ in $[n]$, 
\begin{equation*}
\sett{x_{V_{\lambda_1}}}{x_{V_{\lambda_n}}}\subset \overline{B(x_{V_{\lambda_i}},2s_x/3)}\subset B(x,s_x)
\end{equation*}

For $i$ in $[n]$ let $r_i$ be the minimum value for which $\sett{x_{V_{\lambda_1}}}{x_{V_{\lambda_n}}} \subset \overline{B(x_{V_{\lambda_i}},r_i)}$. We conclude from Proposition \ref{P:close}
that $f(y)$ is in $\overline{B(x_{V_{\lambda_i}},r_i)}\subset B(x,s_x)$. Thus we have shown that the definition of $f(y)$ is independent of the choice of $B(x,s_x)$.

If $y'$ is in $N(y)$ it follows from the definition of $f$ and Corollary \ref{C:PtoP'} that $f(y')$ is the barycenter of $\sett{x_{V_{\lambda_1}}}{x_{V_{\lambda_n}}}$ relative to $\sett{u_{V_{\lambda_1}}(y')}{u_{V_{\lambda_n}}(y')}$. By Theorem \ref{T:con}
 (Part b) and the fact that the $u_V$, the members of the partition of unity, are continuous, we see, by essentially the same argument as in Theorem \ref{T:ar}, that $f$ is continuous at $y$. 
 Now suppose that $x$ is in $\bd(X)$. Let $\epsilon>0$ but require $\epsilon$ to be less than $s_x$. From the definition of $f$, it is clear that there exists a $\delta_1>0$ such that if $y$ is in $\ball{x}{\delta_1}$ then $f(y)$ is in $\ball{x}{s_x}$. Choose $\delta=\min\{\delta_1,\epsilon'/6\}$ where $0<\epsilon'<\epsilon$. By Lemma \ref{L:howfar} and Proposition \ref{P:close}, if $y \in \ball{x}{\delta}$ then $f(y) \in \oball{x}{\epsilon'}\subset \ball{x}{\epsilon}$. Thus $f$ is continuous at $x$ and it follows that $f$ is a retraction of $Y$ onto $X$. By Theorem \ref{T:aranr}, $X$ is an absolute neighborhood retract.  
\end{proof}.

\begin{corollary} \label{T:anrc} Let $(X,d)$ be a $\cat$ space with $k > 0$. Then $X$ is an absolute neighborhood retract.

\end{corollary}

The results above imply certain facts about extensions of maps into $\cat$ spaces and spaces of curvature $\le k$. They also provide information about the fixed point theory of self maps of $\cat$ spaces and spaces of curvature $\le k$. See for example (\cite{bor},\cite{br},\cite{gd}).

Next we consider multifunctions. Let $X$ and $Y$ be metric spaces. A multifunction $\phi$ from $X$ to $Y$, denoted
by $\phi: X \multimap Y$, assigns to each $x$ in $X$ a subset $\phi(x)$ of $Y$. A selection is a function $\fn{f}{X}{Y}$
such that $f(x)$ is in $\phi(x)$ for each $x$. A multifunction $\phi: X \multimap Y$ is lower semi continuous (l.s.c.)
provided that for every open set $U$ in $Y$, the set $ \{ x\vert \phi(x) \cap U \neq \emptyset \}$ is open in $X$. We state the well known theorem of  Michael in the metric space setting (\cite{m1},\cite{m2}, \cite{m3}. \cite{m4},\cite{g}).

\begin{theorem}\label{T:m} Let $X$ be a metric space and $Y$ be a Banach space. Then a l.s.c. multifunction $\phi: X \multimap Y$ with closed convex values has a continuous selection.
\end{theorem}

We formulate a similar theorem in the $\catt$ setting.

\begin{theorem}\label{T:m2} Let $X$ be a metric space and $Y$ be a $\catt$ space. Then a l.s.c. multifunction $\phi: X \multimap Y$ with closed convex values has a continuous selection.
\end{theorem}

Of course convex here is in the geodesic sense. Applying barycenters in place of convex combinations, one finds that the
proof of this theorem is essentially the same as the proof of Theorem \ref{T:m} as given in \cite{g,m1}. This result also follows from properties of geodesic and convex structures as given in \cite{m5}.

Finally we apply barycenters in constructing continuous functions defined on certain quotient spaces of $\catt$ spaces.
Suppose that $X$ is a complete $\catt$ space and $\Gamma$ is a group of isometries of $X$ whose orbits are all bounded. Let $\fn{f}{X}{\bld}$ be a positive valued, bounded and uniformly continuous function. Let $\Gamma/X$ be the quotient space of orbits. For each orbit $O(x)$ in $\Gamma/X$, define $F(O(x))=O(y)$ where y is the barycenter of the support
of $O(x)$ in $X$ relative to $f$. This defines a function $\fn{F}{\Gamma/X}{\Gamma/X}$ and by Theorem \ref{T:con},
it is continuous.

\end{document}